
\magnification = \magstep1

\def\s{\succeq}
\def\p{\preceq}
\def\LP{\left(}
\def\RP{\right)}
\def\LB{\left[}
\def\RB{\right]}

If $f = \sum_n f_n x^n$ and $g = \sum_n g_n x^n$ are power series, we
write $f\s g$ if $f_n\ge g_n$ for all $n$.
Thus if $f_1\s f_2\s 0$ and $g_1\s g_2\s 0$, then $f_1g_1\s f_2g_2$.
If $f\s g$ and $h\s 0$,
then $h(f)\ge h(g)$.  In particular, this is true when $h = \exp$.

\smallskip\noindent
{\it Lemma. If $k > 1$ and $f\s 0$ satisfies $f_m = O(m^{-k})$, then
$\exp(f)_m = O(m^{-k})$.}
\smallskip\noindent
{\bf Proof.} It suffices to prove the result when
$f_m = C m^{-k}$.  Suppose $a(x) = \sum_{m=1}^\infty a_m x^m$ and
$b(x) = \sum_{m=1}^\infty b_m x^m$ are $\s 0$ and
$c(x) = a(x)b(x) = \sum_m c_m x^m$.
If for all $m$,
$a_m \le A m^{-k}$ and $b_m \le B m^{-k}$, then
$$\eqalign{c_m =
\sum_{i=1}^{[m/2]} a_i b_{m-i} + \sum_{i=[m/2]+1}^{m-1} a_i b_{m-i}
&\le B (m/2)^{-k}\sum_{i=1}^\infty A i^{-k} +
A (m/2)^{-k} \sum_{i=1}^\infty B i^{-k} \cr&= 2^{k+1}AB\zeta(k) m^{-k}.}$$
By induction on $n$,
$$(f^n)_m \le C^n 2^{(k+1)(n-1)}\zeta(k)^{n-1} m^{-k},$$
so for $m > 0$,
$$\exp(f)_m \le {\exp\LP 2^{k+1} C\zeta(k)\RP\over 2^{k+1}\zeta(k)} m^{-k}.$$

\smallskip\noindent
{\it Remark.} The theorem does not hold for $k=1$, but one does have the
following variant: if $f_m = o(m^{-1})$ then $\exp(f)_m = o(m^{\epsilon-1})$
for all $\epsilon > 0$.
\smallskip
Now we apply this lemma to the problem of partitions.  Fix an integer
$n\ge 2$.  Let $p_k$ denote the probability that a random permutation
on $k$ letters, drawn from a the uniform distribution on $S_k$, has
an $n^{\rm th}$ root.  The condition is equivalent to the statement
that the number of $r$-cycles in the permutation is divisible by
$d_n(r)$, where $d = d_n(r)$ denotes the largest divisor of $n$ such
that $r$ and $n\over d$ are relatively prime.
By a standard generating function argument
$$p(x) = \sum_{k=0}^\infty p_k x^k = \prod_{r=1}^\infty
\exp_{d_n(r)}(x^r/r),$$
where
$$\exp_k(x) := \sum_{i=0}^\infty {x^{ik}\over (ik)!}\p
\exp(x^k/k!).$$
Therefore,
$$p(x)\p q_1(x)q_2(x) := \exp\LP\sum_{(r,n)=1}{x^r\over r}\RP
\exp\LP\sum_{(r,n)> 1}x^{r d_n(r)}\over r^2\RP.$$
We can break up the logarithm of the second multiplicand, $\log(q_2(x))$,
into a finite sum
over residue classes of $n$; each such sum has $O(m^{-2})$
coefficients, so $\log(q_2)$ and hence $q_2$ has $O(m^{-2})$ coefficients.
The $q_1(x)$ term can
be expressed by the M\"obius inversion formula as
$$\prod_{d\vert n} (1 - x^d)^{-\mu(d)/d}$$
%
By Darboux's lemma [Kn-W], the coefficients of this power series are
asymptotic to a multiple of $m^{\phi(n)-n\over n}$.

Finally, we observe that if $0<\alpha<1$, $a(x)\s 0$ is any power series
with coefficients $a_m \sim C m^{-\alpha}$, and $b(x)\s 0$
is a power series with coefficients $b_m = O(m^{-2})$, then
the coefficients of $c(x) = a(x)b(x)$ are aymptotic to
$(C\sum_i b_i) m^{-\alpha}$.
As the $b_i$ are non-negative, $a(x)b(x)\s a(x)(b_0+b_1x+\cdots+b_kx^k)$.
It follows that
$$\liminf_{m\to\infty} {c_m\over m^\alpha}\ge 
\lim_{m\to\infty} C\LP b_0 + b_1(1-1/m)^{-\alpha} + \cdots +
b_k(1-k/m)^{-\alpha}\RP =
C\sum_{i=0}^k b_i,$$
and sending $k\to\infty$
$$\liminf_{m\to\infty} {c_m\over m^\alpha}\ge C\sum_{i=0}^\infty b_i.$$
On the other hand, fixing $\beta\in(\alpha,1)$,
$$c_n = \sum_{i=0}^{\LB n^\beta\RB} b_i a_{n-i} +
\sum_{i=\LB n^\beta\RB+1}^n
b_i a_{n-i} < \LP\LP\sum_{i=0}^\infty b_i\RP\sup_{j\in(n-n^\beta,n)} a_j\RP
+ \LP\sum_{i=\LB n^\beta\RB+1}^\infty b_i\RP\sup_j a_j.$$
As $a_j$ is absolutely bounded, 
$$\sum_{i=[n^\beta]+1}^\infty b_i = O(n^{-\beta}),$$
and $a_m\sim C m^{-\alpha}$,
$$\limsup_{n\to\infty} {c_n\over n^\alpha}\le C\sum_{i=0}^\infty b_i.$$
We conclude that $p_m\sim C m^{\phi(n)-n\over n}$ for any fixed $m$.

\end